 \documentclass[12pt]{article}
 \usepackage{graphicx}

 \usepackage{amssymb}

 \font\smallbold=cmbx9 
  \font\smallbf=cmbx10  scaled \magstep1
 \font\bigtext=cmcsc10 scaled\magstep1
 \font \smallrm=cmr10

 \def \t{\theta}
 \def \T{\Theta}

 \def \e{\varepsilon}
 \def \wt{\widetilde}
\def \dist{\rm dist}

\begin{document}
  \thispagestyle{plain}
  \setlength{\baselineskip}{0.5truecm}
  \vskip40pt
 \centerline
 {\smallbf On $ \Lambda$-positioning of an arc between two parallel support lines}
  \vskip8pt \centerline{\smallbold Y. Movshovich}
 \begin{abstract} {\smallrm
 We show that  a rectifiable plane arc $g $
 has two parallel support lines and a triple of consecutive points
 $g(r),\ g(s),\ g(t),\ \ r<s<t, $ so that $g(s)$ lies on one
 line, while $g(r)$ and $ g(t) $ lie on the other.
  If the arc is simple, such a pair of lines is unique.}
  \end{abstract}
 \baselineskip18pt
   \vskip 1 pt  \textbf{Introduction.} In the articles \cite{CM, M} 
      we had to use the result  in \cite{MWW}: 
 {\it If a convex set  covers any simple polygonal unit arc, it covers any unit arc}. In \cite{CM, M} 
 the requirement on an arc to be simple and polygonal was used only 
   in Theorem 5.1 of \cite{CM} 
    establishing that {\it any simple polygonal arc assumes
   a so-called $\Lambda$-configuration} (Figure 1). 
   Two proofs of Theorem 5.1 exist for simple arcs: one by Y. M. (Geometry Seminar, UIUC, 2009) and
   the other by  R. Alexander, J. E. Wetzel, W. Wichiramala in their recently submitted paper "The $\Lambda$-property of a simple arc".
      In this note we prove Theorem 5.1 of \cite{CM} omitting both requirements: simple and polygonal. 
 \vskip2pt
 Given a parametrization $g(s),\ s\in[0,1],$ for points
   $p=g(s_1),\ q=g(s_2)$ with  $s_1<s_2,$ we say that $p$ {\it precedes} $q$ and write it as $p \prec q.$
  Points $p_1,\ p_2,\ p_3 $ form
 a {\it triple of consecutive points} if either
 $$p_1\prec p_2\prec p_3 \quad {\rm or} \quad p_3\prec p_2\prec p_1. \leqno (1) $$
  We are seeking a pair of parallel support lines with and a triple  of consecutive points $p_1,\ p_2,\ p_3 $  such that
   $p_1,  p_3$ lie on one line and  $p_2$ lies on the other.  
  \begin{center}  \includegraphics[scale=0.4]{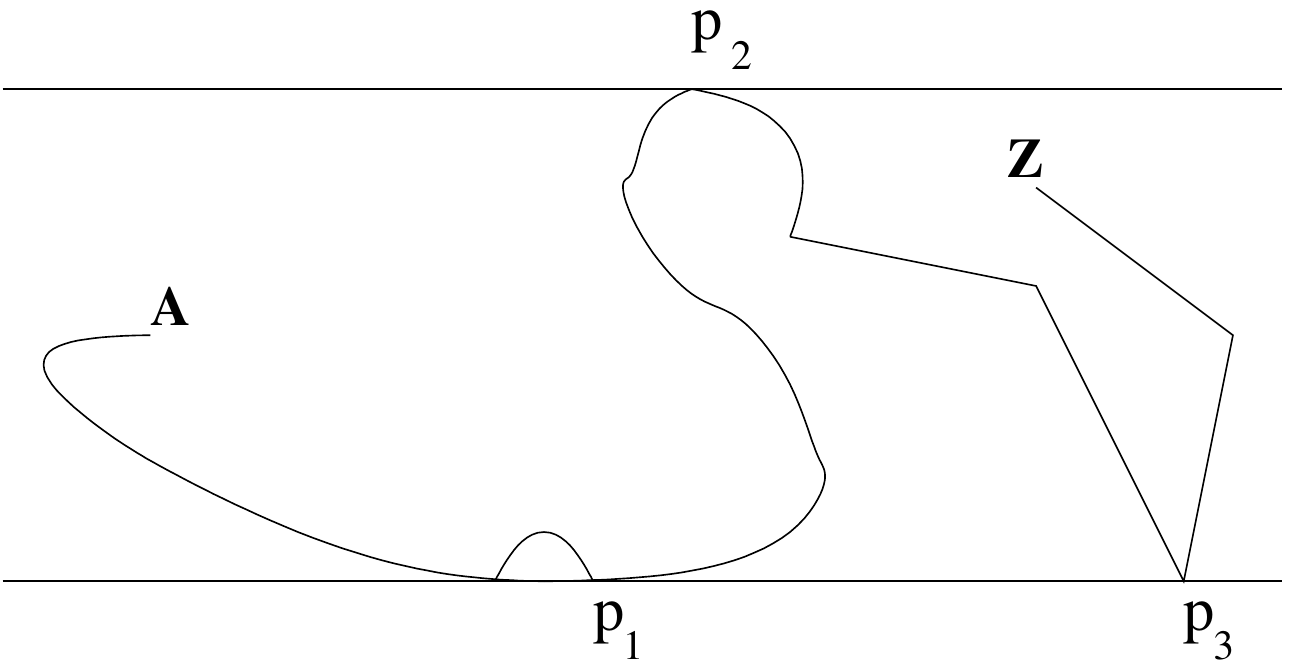}
  \par Figure 1. \quad  {\bf \boldmath   $  \Lambda$-configuration}.
  \end{center}
This is clearly true when $g$ is closed or is a straight segment. 
\vskip2pt
 
  %
 \par {\bf Lines $l(\t), $ walls $W(\t)$ and touch points in Figure 2.}
 Let  $g(s) :[0,1]\to R^2$ be an open rectifiable plane arc in the horizontal $(x,y)$-plane, whose  thickness $h$  is positive. 
 We denote by $l(\t)$ a counter-clockwise oriented  support line for  $g$  in the direction of polar angle $\t$  
   and assume  that $g$ lies above $l(0).$
 Let $G(s)=\big(g(s),s\big):[0,1]\to R^2\times[0,1]$ be a simple lift of $g$ to $R^3$ with  $G(0)=A,$  $G(1)=Z.$
 Given a direction $\t,$ let us denote by $W(\t)$ the support plane (wall)
to $G$ through $l(\t)$ orthogonal to the $(x,y)$-plane.
The touch points $G\cap W(\t)$ are denoted by $I_k.$
 Note that a set of touch points of a wall $ W(\t)$ is its compact subset.
 The wall $W(\t)$ has  its lowest most left $I_1(\t)$ and highest most right $C(\t)$ touch points with respect to the counter-clockwise orientation of the support line $l(\t).$ 
 Let $v(\t)$ be the vertical unit segment  through  $C(\t).$ 
  \begin{center}  \includegraphics[scale=0.44]{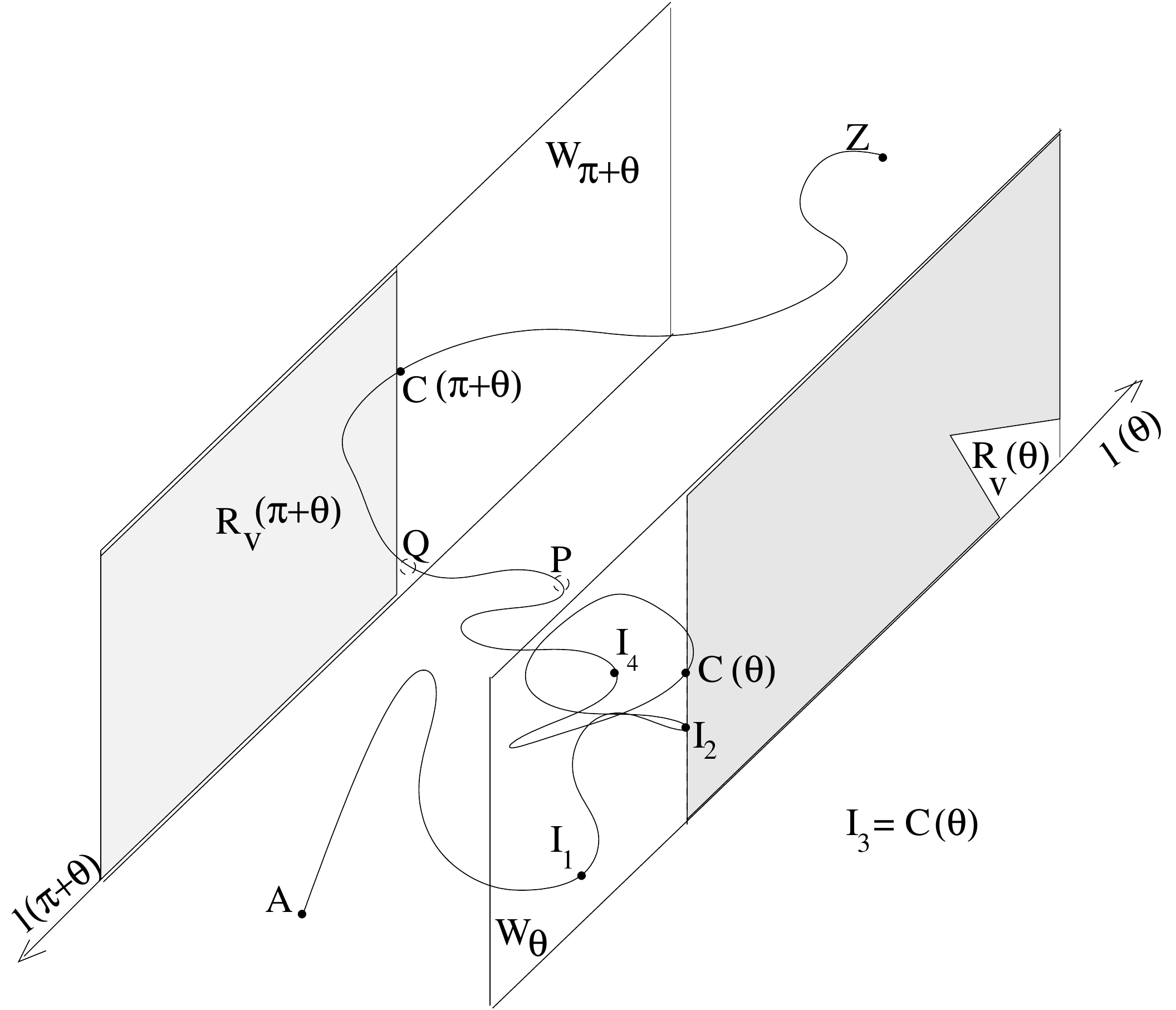}
  \par {\boldmath \smallbold Touch points $\scriptstyle I_k(\t)\, \in \,W(\t)$,     
  $\scriptstyle I_m(\t+\pi)= C(\t+\pi)\, \in \, W(\t+\pi)$. 
  For all $\scriptstyle k,m,$ $\scriptstyle I_k(\t)\prec  I_m(\t+\pi)$.}
   \par Figure 2
  \end{center}

 \vskip2pt
  {\bf Rotation of a support wall as a map.}
 Note that for each $\t\in [0,2\pi),$ there is one  support line $l(\t)$ and thus one support wall $W(\t).$ 
  Referring to a rotation of a support line around  $g$ as $\t$ changes from $0$ to $2\pi,$
 we think of two maps from $ [0,2\pi):$ one is the map to the set of all lines in  the $(x,y)$-plane and the other is 
 the map to all planes in the $(x,y,z)$-space.
 The images of  such maps are the set of oriented support lines $l(\t)$ of $g$ and
 the set of support walls  $W(\t)$ through $l(\t).$ 
  Geometrically the first map is represented by  a line  moving in a plane so that it
 coincides with $l(\t)$ for each $\t,$ while the second map is represented by
  a plane moving in the space so that it coincides with $W(\t)$ for each $\t.$
    \vskip 6 pt
    {\bf Local stability of a non-$\Lambda$-configuration.}  
    A subarc of $g$ between points $X$ and $Y$ is denoted by $\widetilde{XY}.$
    \vskip6pt
    \textbf{Lemma 1.}
    {\it  The following set
    $ \displaystyle \T=\{\t\ge0:\ C(\t)\prec C(\t+\pi) \ \}$ is a half-open interval.
  That is
  if \ $\t\in \T,$ \ then there exist \ $\delta$ \ so that \
 $\t+\e\in \T$ \ for any \ $\e<\delta.$ }
 \vskip1pt
  {\it Proof.} To keep our proof transparent, we assume that $v(\t)$ has only finitely many touch points and we will use a particular configuration of Figure 2.  
  Let points $P,Q\in \widetilde{I_4D}$ be so that
  $$I_4(\t)\prec P\prec Q\prec C(\t+\pi) \ \ \  and\ \ \
 length(\widetilde{I_4P})=length(\widetilde{QD})={1\over3}h.
 \ \ \ \ \  \leqno (2)$$
 Denote by $R_v(\t)\subset W(\t)$ a half-plane of points to the right of $v(\t).$
  Let
 $$\sigma(\t) = \min\Big\{  \dist \big[\widetilde{PZ} , R_v(\t) \big],\
 \dist \big[ \widetilde{AQ} , R_v(\pi+\t) \big] \Big\}\ . $$
 Then $\sigma>0$ and  $ \dist \big(\widetilde{PZ} ,v(\t) \big)$ and
 $ \dist \big(\widetilde{AQ} ,v(\t+\pi) \big)$ are $\ge\sigma.$
 We take any
 $$\displaystyle \e\ <\ \delta={\sigma(\t)\over88\ diameter(G)}. $$
 The  obstacles to
 the counter-clockwise rotation of the walls by $\e$ could be
 only subarcs  $\widetilde{QZ}$ and $\widetilde{AP}.$
  Therefore, 
  $ C(\t+\e)\in\widetilde{AP}$ while $ D(\t+\e)\in\widetilde{QZ}$ and hence
 by (2),\ 
 $C(\t+\e)\prec D(\t+\e).$
 \ $  \blacklozenge $
 %
  \vskip 9 pt
 \textbf{Theorem 1.}
  {\it
   Let $g$ be an open rectifiable arc with a thickness $h>0.$
 Then there exist support lines $l(\t)$ and $l(\pi+\t)$
 containing a triple of consecutive  points
 $p_1,\ p_2,\ p_3 $ of $g$ with the lone middle point $p_2.$
  }
   \vskip1pt {\it Proof.} 
   We may assume that
   $A\prec C(0),C(\pi)\prec Z.$ (Figure 2).
  If a triple of the theorem  exists in the strip between $l(0)$ and $l(\pi)$  then   $C(\pi)\prec C(0)$ (Figure 1).
   Otherwise, $$C(0)\prec C(\pi) . \leqno (3)$$ 
   By Lemma 1 this property is  locally  stable and so if
  $\T$ is the set given by this lemma  and  {\boldmath $\t_T=lub(\T)$},
 then   \ $\t_T\notin \T .$ That is   $C(\t_T+\pi)\prec C(\t_T).$
 However,  limits of most right touch points in $\T$ preserve this property of
 $\T:$
 $$\lim_{\t\uparrow\t_T}C(\t_T)\prec\lim_{\t\uparrow\t_T}C(\t_T+\pi), \leqno (4)$$
 because they are separated in distance by $h.$ Indeed,
 \vskip2pt {\bf Small rotations around $G$ by a small angle $\psi$.} Configurations of local behavior near a touch point
 are given in Figure 3:  
  \begin{center}
  \includegraphics[scale=0.5]{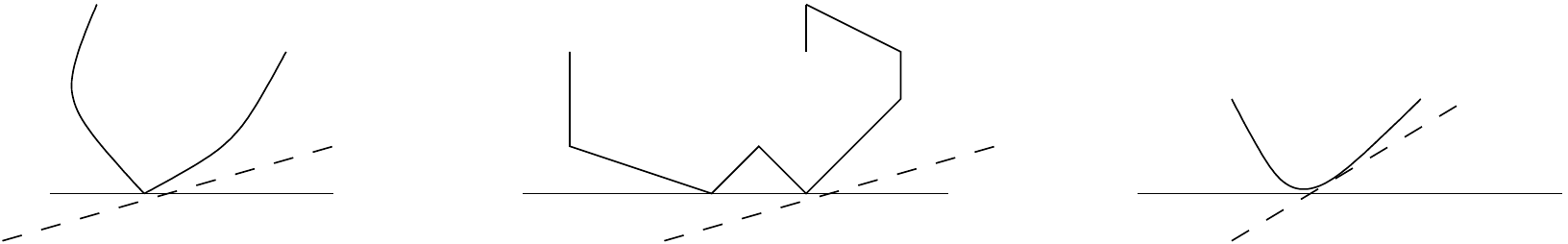}  
   \par {\boldmath \smallbold $\scriptstyle I(\t+\psi)=I(\t)$ in rotations without an obstacle,      
 $\scriptstyle \dist\big(I(\t+\psi), I(\t)\big)$ is small in rotations around smooth convex arcs.}
   \par Figure 3
  \end{center}

 \vskip2pt
 Thus one or both  $C(\t_T)$ or  $C(\t_T+\pi)$ are not  limits in (4).
  Suppose that $C(\t_T)$  
   is not equal to
    $\displaystyle  \lim_{\t\uparrow\t_T}C(\t_T)$ 
  Then the triples
   $p_1\prec p_2\prec p_3$ satisfying the theorem are either
 $$\matrix{ \displaystyle \lim_{\t\uparrow\t_T}C(\t_T)\prec\lim_{\t\uparrow\t_T}C(\t_T+\pi)\prec C(\t_T) \cr
   \cr
\displaystyle   C(\t_T+\pi)\prec C(\t_T)\prec\lim_{\t\uparrow\t_T}C(\t_T+\pi).
 } \leqno\hbox{or} $$
   {\bf \boldmath An assumption that there were no  $\t_T\le\pi$ leads to a contradiction.} 
   Suppose that $\t_T >\pi,$ that is  $C(\pi)\prec C(\pi+\pi).$
    On the other hand, after a rotation by $\pi,$
   we  arrive to the initial position of $g$ between the same parallel support lines.
   In this configuration, the highest most right touch point on the wall $W(\pi)$ 
   is a successor to such point on  the wall $W(\pi+\pi)=W(0)$ and  (3) is true:
   $ C(\pi+\pi)=C(0) \, \prec \, C(\pi) .$ 
   That is a contradiction.
      $\blacksquare$
%
  \vskip6pt 
  {\bf Corollary} (W. Wichirimala).
   {\it If $g$ is simple, then the pair of support lines with $g$ positioned as in Figure 4 is unique.}
  \begin{center}
  \includegraphics[scale=0.44]{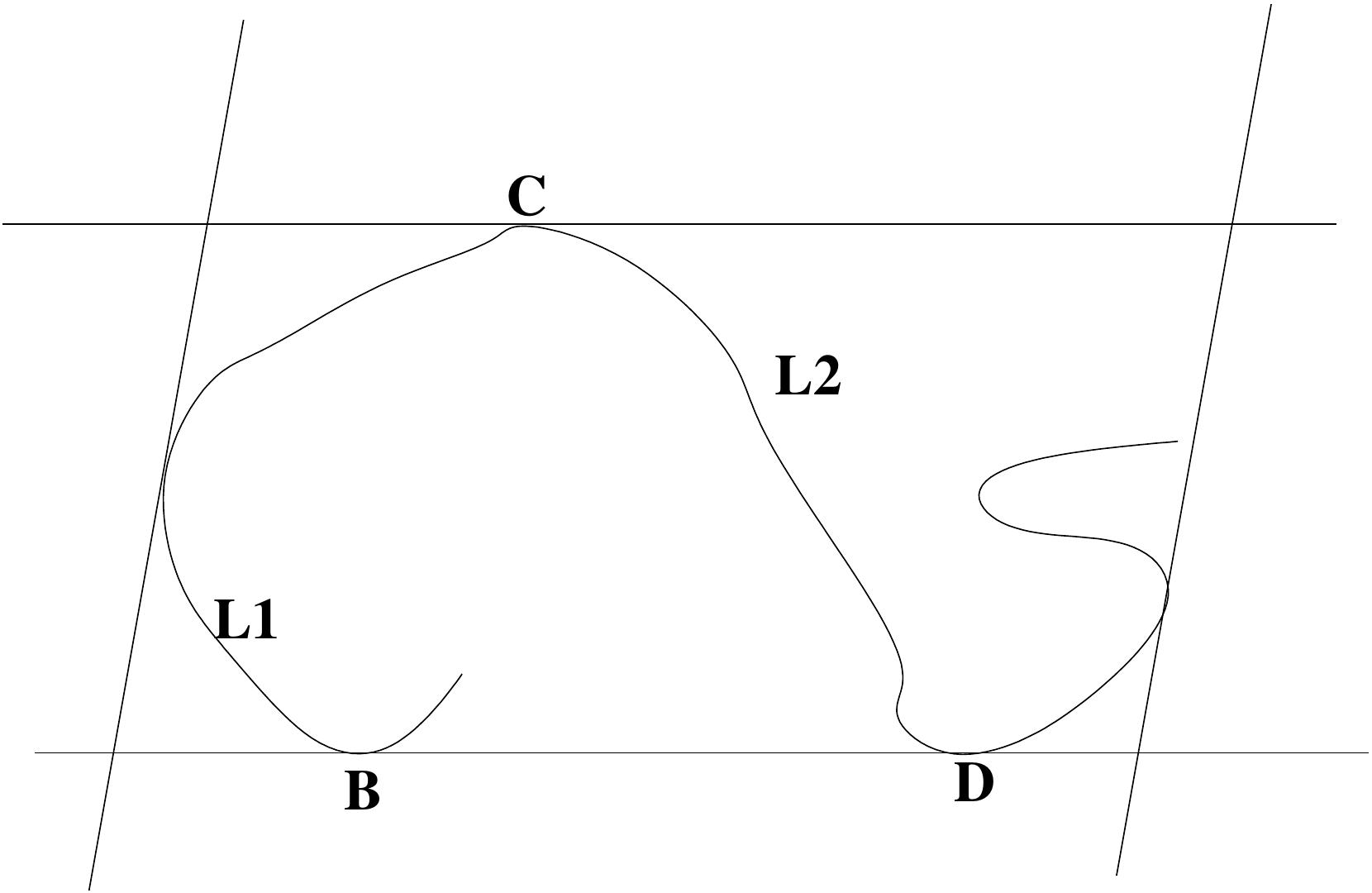}
   \par Figure 4
  \end{center}
{\it Proof.}
 In Figure 4, the point $C$  divides the curve into two curves $L1$ and $L2$
 where each point on $L1$ is parametrically precedes each point on $L2.$ 
 Consider  a different pair of two parallel support lines. All touch points on one, say, the left line belong to $L1,$ therefore
 none of them can be parametrically between two touch points on the right line
 belonging to $L2.$ 
   \ $  \blacklozenge $

                  \end{document}